\documentclass[12pt,preprint,nofootinbib,nobibnotes,nonumbib]{article}
\usepackage{amssymb}
\usepackage{amsmath}
\usepackage{hyperref}
\usepackage{amscd}
\usepackage{graphicx}

\begin{document}

\title{${\mathbb Z}_n$-equivariant $K$-theory}

\author{Larry B. Schweitzer} 

\date{May 7, 2015}

\maketitle

\begin{abstract}
We construct a sequence of $n-1$ cyclic exact sequences that 
can be used to compute the $K$-theory of 
the $C^\star$-algebra crossed product $A \ltimes {\mathbb Z}_n$.
\end{abstract}

\eject

\tableofcontents

\section{Introduction}
W. L. Paschke's exact sequence
for computing the $K$-theory of the $C^\star$-algebra crossed
product $A \ltimes_\alpha {\mathbb Z}_2$  is 
\begin{equation}
\begin{CD}
K_0(A_0) @>>> K_0(A \ltimes_\alpha {\mathbb Z}_2) @>>> K_0(A_0/J) \\
@AAA & & @VVV\\
K_1(A_0/J)  @<<<  K_1(A \ltimes_\alpha {\mathbb Z}_2) @<<< K_1(A_0)
\end{CD}, \label{eq:Paschke}
\end{equation}
where $\alpha$ is an action of ${\mathbb Z}_2$ on the $C^\star$-aglebra
$A$, $A_0$ is the fixed point algebra of $\alpha$,
and  $J$ is the ideal of $A_0$ generated by
$\{\, xy \, | \, x, y \in A_1 \, \}$, where
$A_1= \{ \, x \in A \,|\, \alpha(x) = -x \, \}$ [Pas, 1985].
See [Ich, 1989] for a generalization of (\ref{eq:Paschke})
to the case of crossed products by
compact groups.
We present a generalization of (\ref{eq:Paschke})
to crossed products by
${\mathbb Z}_n$, where $n$ is any positive
integer, $n \in \{2,3,4, \dots\}$.

\section{Embedding the crossed product as a subalgebra of $M_n(A)$}
Let $\alpha$ be an action of ${\mathbb Z}_n$ on a $C^\star$-algebra $A$.
Let $\xi$ be a primitive $n$th root of unity.  For $0 \leq k \leq n-1$,
the $k$th eigenspace of $\alpha$ is
\begin{equation}
A_k = \{\, x \in A \, | \, \alpha(x)=\xi^k \alpha(x)\, \}.
\label{eq:Asubk}
\end{equation}
Then $A_k^\star = A_{n-k}$.  Note that $A_0$ is the fixed point algebra
of $\alpha$, and each $A_k$ is an $A_0$-bimodule.
Denote by $A_j A_k$ the closed linear span of
$\{\, xy\, |\, x \in A_j, y \in A_k  \, \}$.
Note that $A_j A_k \subseteq A_{j+k}$ and $A_k A_{n-k}$
is a two-sided ideal in $A_0$.

Let $P_k \colon A \rightarrow A_k$ be the projection
\begin{equation}
P_k(x) = {1 \over n} \sum_{i=0}^{n-1} \xi^{-ki} \alpha^i(x), \qquad x \in A.
\end{equation}
Then $\alpha(P_k(x)) = \xi^k P_k(x)$, $P_k$ is onto $A_k$,
$P_j P_k = \delta_{jk} P_k$, and $\sum_{k=0}^{n-1} P_k = I$.

Elements of the crossed product $A \ltimes_\alpha {\mathbb Z}_n$
are denoted by $\sum_{i=0}^{n-1} a_i \lambda_i$, $a_i \in A$,
where $i \mapsto \lambda_i$ a unitary representation of ${\mathbb Z}_n$.
Multiplication is given by
\begin{equation}
\bigl( a_i \lambda_i \bigr) \bigl( b_j \lambda_j \bigr)
= a \alpha^i(b_j) \lambda_{i+j}.
\end{equation}

\vskip\baselineskip
\noindent{\bf Lemma 1.}\ 
{\it The linear map $\theta \colon
A \ltimes_\alpha {\mathbb Z}_n \rightarrow M_n(A)$ defined by
\begin{equation}
\theta(a_i \lambda_i) =
\begin{pmatrix} P_0(a_i) & \xi^{-i} P_1(a_i) & \dots
& \xi^{-i(n-1)}P_{n-1}(a_i) \\
P_{n-1}(a_i) & \xi^{-i} P_0(a_i) & \dots
& \xi^{-i(n-1)}P_{n-2}(a_i) \\
\vdots & \vdots & \ddots & \vdots \\
P_1(a_i) & \xi^{-i} P_2(a_i) & \dots & \xi^{-i(n-1)}P_0(a_i) \end{pmatrix}
\end{equation}
is a $\star$-isomorphism of
$A \ltimes_\alpha {\mathbb Z}_n$ with the $C^\star$-subalgebra
\begin{equation}
\begin{pmatrix} A_0 & A_1 & \dots
& A_{n-1} \\
A_{n-1} & A_0 & \dots
& A_{n-2} \\
\vdots & \vdots & \ddots & \vdots \\
A_1 & A_2 & \dots & A_0 \end{pmatrix}
\label{eq:theSubalg}
\end{equation}
of $M_n(A)$.} 

\noindent{\bf Proof:}
First show $\theta$ is an algebra homomorphism.
Note $\bigl( \theta(a_i \lambda_i) \bigr)_{kl}
= \xi^{-il} P_{l-k}(a_i)$, $0 \leq k, l, \leq n-1$.
Thus
\begin{eqnarray}
\bigl( \theta(a_i \lambda_i) \theta(b_j\lambda_j) \bigr)_{kl}
&=& \sum_{m=0}^{n-1}
\bigl( \theta(a_i \lambda_i) \bigr)_{km}
\bigl( \theta(b_j \lambda_j) \bigr)_{ml} \nonumber \\
&=& \sum_{m=0}^{n-1}
\xi^{-im} P_{m-k}(a_i)
\xi^{-jl} P_{l-m}(b_j) \nonumber \\
&=& \sum_{m=0}^{n-1}
\xi^{-im - kl}
\biggl( {1 \over n} \sum_{p=0}^{n-1} \xi^{-(m-k)p} \alpha^p(a_i)\biggr)
\biggl( {1 \over n} \sum_{q=0}^{n-1} \xi^{-(l-m)q} \alpha^q(b_j)\biggr)
\nonumber \\
&=&
 {1 \over n^2}
\sum_{m,p,q=0}^{n-1}
\xi^{-im - jl -(m-k)p -(l-m)q}
\alpha^p(a_i)
 \alpha^q(b_j)
\nonumber \\
&=&
 {1 \over n^2}
\sum_{p,q=0}^{n-1}
\xi^{- jl +kp -lq}
\sum_{m=0}^{n-1}
\biggl( \xi^{m(- i- p + q)}  \biggr)
\alpha^p(a_i)
 \alpha^q(b_j)
\nonumber \\
&=&
 {1 \over n^2}
\sum_{p,q=0}^{n-1}
\xi^{- jl +kp -lq}
\biggl(
n \delta_{- i- p + q,0}  \biggr)
\alpha^p(a_i)
 \alpha^q(b_j)
\nonumber \\
&=&
 {1 \over n}
\sum_{p=0}^{n-1}
\xi^{- jl +kp -l(i+p)}
\alpha^p(a_i)
 \alpha^{i+p}(b_j)
\qquad \mbox{$-i-p+q=0\Rightarrow q=i+p$}
\nonumber \\
&=&
\xi^{-(i+j)l}
\biggl(
 {1 \over n}
\sum_{p=0}^{n-1}
\xi^{-(l-k)p}
\alpha^p\bigl(a_i \alpha^i(b_j)\bigr)
\biggr)
\nonumber \\
&=&
\xi^{-(i+j)l}
P_{l-k}\biggl(a_i \alpha^i(b_j)\biggr)
\qquad \qquad\qquad \mbox{definition of $P_{l-k}$}
\nonumber \\
&=&
\theta \biggl( a_i \alpha^i(b_j) \lambda_{i+j} \biggr)_{kl}
\nonumber \\
&=&
\theta \biggl( \bigl(a_i \lambda_i \bigr)
\bigl(b_j \lambda_j \bigr) \biggr)_{kl}.
\end{eqnarray}
It follows that
$ \theta \biggl( \bigl(\sum a_i \lambda_i \bigr)
\bigl(\sum b_j \lambda_j \bigr) \biggr)
=
\theta\bigl(\sum a_i \lambda_i\bigr)
\theta\bigl(\sum b_j\lambda_j\bigr)
$.
Similarly, one checks that
$
\theta \biggl( \bigl(\sum a_i \lambda_i \bigr)^\star
 \biggr)
=
\biggl(\theta\bigl(\sum a_i \lambda_i\bigr)
 \biggr)^\star$.

To see that $\theta$ is one-to-one, note that an
inverse map exists.  This is because for
$\sum a_j\lambda_j \in A \ltimes_\alpha {\Bbb Z}_n$,
we can retrieve each $a_j$ via the formula
\begin{equation}
{1\over n} \sum_{k,l=0}^{n-1} \xi^{li}
\bigl( \theta \bigl( \sum a_j \lambda_j \bigr) \bigr)_{kl} = a_i.
\label{eq:inverse}
\end{equation}
The map $\theta$ is onto since this inverse map
is well defined if we
replace $\theta\bigl(\sum a_j \lambda_j\bigr)$ by an
arbitrary element of the subalgebra (\ref{eq:theSubalg})
of $M_n(A)$.
$\square$

\vskip\baselineskip
Let $B_0$ denote the subalgebra (\ref{eq:theSubalg}) of Lemma 1.
More generally, let
$B_k$ be the $k$th lower right corner of $B_0$,
\begin{equation}
B_k=
\begin{pmatrix}
A_0 & A_1 & \dots & A_{n-(k+1)} \\
A_{n-1} & A_0 & \dots & A_{n-(k+2)} \\
\vdots & \vdots & \ddots & \vdots \\
A_{k+1} & A_{k+2} & \dots & A_0
\end{pmatrix},
\end{equation}
for $0 \leq k \leq n-1$, which
is $n-k$ by $n-k$.
Note $B_{n-1} = A_0$.  For $0 \leq k \leq n-2$,
let $I_k$ be
the two-sided closed ideal
$\sum_{i=1}^{n-1-k} A_i A_{n-i}$
in $A_0$.
Note that
for $0\leq k \leq n-2$,
\begin{equation}
J_k=
\begin{pmatrix}
I_k & A_1 & \dots & A_{n-(k+1)} \\
A_{n-1} & A_0 & \dots & A_{n-(k+2)} \\
\vdots & \vdots & \ddots & \vdots \\
A_{k+1} & A_{k+2} & \dots & A_0
\end{pmatrix}
\end{equation}
is a closed two-sided ideal of the $C^\star$-algebra
$B_k$.
Thus for $0\leq k \leq n-2$, we have $n-1$ exact sequences
\begin{equation}
0 \;\xrightarrow\; J_k \;\xrightarrow\; B_k \;\xrightarrow\; B_k/J_k \cong
A_0/I_k\; \xrightarrow\; 0.
\label{eq:algExacts}
\end{equation}

\section{Exact sequences for the $K$-theory}

From the $n-1$ exact sequences (\ref{eq:algExacts})
of $C^\star$-algebras,
we get $n-1$ cyclic exact sequences in $K$-theory 
\begin{equation}
\begin{CD}
K_0(J_k) @>>> K_0(B_k) @>>> K_0(A_0/I_k) \\
@AAA & & @VVV\\
K_1(A_0/I_k)  @<<<  K_1(B_k) @<<< K_1(J_k)
\end{CD},
\label{eq:nm1Exact}
\end{equation}
for $0\leq k \leq n-2$.

\vskip\baselineskip
\noindent{\bf Lemma 2.}\ 
{\it Assume $A$, $I_0$, $I_1$, $\dots$ $I_{n-2}$ all
have strictly positive elements.
Then $A_0$ has a countable approximate unit for $A$.
For each $0 \leq k \leq n-2$, $J_k$
has a countable approximate unit, $J_k$ has a strictly
positive element, and $B_{k+1}$
is a full corner of $J_k$.}

\noindent{\bf Proof:}
Let $\{e_n\}_{n=0}^\infty$ be an approximate unit for $I_k$,
and $\{f_n\}_{n=0}^\infty$ an approximate unit for $A$ in $A_0$.
Then
\begin{equation}
\begin{pmatrix}
e_n & 0 & \dots & 0 \\
0 & f_n & \dots & 0 \\
\vdots & \vdots & \ddots & 0 \\
0 & 0 & \dots & f_n
\end{pmatrix}
\end{equation}
is an approximate unit for
\begin{equation}
J_k=
\begin{pmatrix}
I_k & A_1 & \dots & A_{n-(k+1)} \\
A_{n-1} & A_0 & \dots & A_{n-(k+2)} \\
\vdots & \vdots & \ddots & \vdots \\
A_{k+1} & A_{k+2} & \dots & A_0
\end{pmatrix}.
\end{equation}
It suffices to show that $e_n a \rightarrow a$ for all
$a \in A_1, \dots, A_{n-(k+1)}$.
Note that
\begin{equation}
\| e_n a - a \|^2 = \| (e_n a - a)(e_n a - a)^\star \| =
\| e_n a a^\star e_n - e_n a a^\star - a a^\star e_n + a a^\star \|.
\label{eq:normEst}
\end{equation}
Now $I_k = A_1 A_{n-1} + \dots A_{(n-1)-k} A_{k+1}
= A_1 A_1^\star + \dots A_{(n-1)-k} A^\star_{(n-1)-k}$
so $aa^\star \in I_k$.  Thus (\ref{eq:normEst}) $\rightarrow 0$
as $n \rightarrow \infty$.

$B_{k+1}$ is a corner of $J_k$ since
\begin{equation}
\begin{pmatrix}
0 & 0 & \dots & 0 \\
0 & 1 & \dots & 0 \\
\vdots & \vdots & \ddots & \vdots \\
0 & 0 & \dots & 1
\end{pmatrix}
J_k
\begin{pmatrix}
0 & 0 & \dots & 0 \\
0 & 1 & \dots & 0 \\
\vdots & \vdots & \ddots & \vdots \\
0 & 0 & \dots & 1
\end{pmatrix}
=
\begin{pmatrix}
0 &  0  \\
0 &  B_{k+1}
\end{pmatrix}.
\end{equation}

$B_{k+1}$ is full since
\begin{eqnarray}
J_k B_{k+1} J_k & \supseteq &
\begin{pmatrix}
0 & A_1 A_0 & \dots & A_{n-(k+1)} A_0 \\
0 & A_0^2 & \dots & A_{n-(k+2)} A_0 \\
\vdots & \vdots & \ddots & \vdots \\
0 & A_{k+2} A_0 & \dots &  A_0^2
\end{pmatrix}
J_k  \\
& \supseteq &
\begin{pmatrix}
I_k & A_1 & \dots & A_{n-(k+1)} \\
A_{n-1} & A_0 & \dots & A_{n-(k+2)} \\
\vdots & \vdots & \ddots & \vdots \\
A_{k+1} & A_{k+2} & \dots & A_0
\end{pmatrix} \\
& = &
J_k.
\end{eqnarray}
Here I used $A_i A_0 = A_i$ for $0 \leq i \leq n-1$.
Let $\{ e_{\gamma} \}$ be an approximate unit for $A$.
Note that $P_0 e_\gamma$ is also an approximate
unit for $A$, since
\begin{equation}
\lim_\gamma (P_0 e_\gamma) a
= \lim_\gamma {1 \over n} \sum \alpha^i(e_\gamma) a
= {1\over n} \sum a = a.
\end{equation}
Thus $A_0$ contains an approximate unit for $A$, so $A_i A_0=A_i$
for each $i$, $0 \leq i \leq n-1$.
$\square$

\vskip\baselineskip
By Lemma 2 and Corollary 2.6 of [Bro, 1977], we have
\begin{equation}
J_k \otimes {\cal K} \cong B_{k+1} \otimes {\cal K},
\end{equation}
for $0 \leq k \leq n-2$, where $\cal K$ is the compact
operators on a separable Hilbert space.
Thus $K_* (J_k) = K_*(B_{k+1})$ and from our
$n-1$ exact sequences (\ref{eq:nm1Exact}),
we get the $n-1$ exact sequences
\begin{equation}
\begin{CD}
K_0(B_1) @>>> K_0(A \ltimes_\alpha {\mathbb Z}_n) @>>> K_0(A_0/I_0) \\
@AAA & & @VVV\\
K_1(A_0/I_0)  @<<<  K_1( A \ltimes_\alpha {\mathbb Z}_n) @<<< K_1(B_1)
\end{CD} \label{eq:fullCrossCD}
\end{equation}
\begin{equation}
\begin{CD}
K_0(B_2) @>>> K_0(B_1) @>>> K_0(A_0/I_1) \\
@AAA & & @VVV\\
K_1(A_0/I_0)  @<<<  K_1(B_1) @<<< K_1(B_2)
\end{CD}
\end{equation}
\begin{equation}
\vdots \nonumber
\end{equation}
\begin{equation}
\begin{CD}
K_0(B_{n-2}) @>>> K_0(B_{n-3}) @>>> K_0(A_0/I_{n-3}) \\
@AAA & & @VVV\\
K_1(A_0/I_{n-3})  @<<<  K_1(B_{n-3}) @<<< K_1(B_{n-2})
\end{CD}
\label{eq:nm3Exact}
\end{equation}
\begin{equation}
\begin{CD}
K_0(A_0) @>>> K_0(B_{n-2}) @>>> K_0(A_0/I_{n-2}) \\
@AAA & & @VVV\\
K_1(A_0/I_{n-2})  @<<<  K_1(B_{n-2}) @<<< K_1(A_0)
\end{CD},
\label{eq:nm2Exact}
\end{equation}
where I've identified $B_0 \cong A \ltimes_\alpha {\mathbb Z}_n$
in (\ref{eq:fullCrossCD})
and $B_{n-1} \cong A_0$ in (\ref{eq:nm2Exact}).
These exact sequences give a recursive procedure for computing
$K_*(A \ltimes_\alpha {\mathbb Z}_n)$ in terms of
$K_*(A_0)$, $K_*(A_0/I_{n-2})$, $\dots$ $K_*(A_0/I_0)$ as follows.
Use (\ref{eq:nm2Exact}) to compute $K_*(B_{n-2})$
in terms of $K_*(A_0)$ and $K_*(A_0/I_{n-2})$.
Then use (\ref{eq:nm3Exact}), $K_*(A/I_{n-3})$,
and $K_*(B_{n-2})$ to compute $K_*(B_{n-3})$, etc.
Finally use (\ref{eq:fullCrossCD}), $K_*(B_1)$, and
$K_*(A_0/I_0)$ to compute
$K_*( A \ltimes_\alpha {\mathbb Z}_n )$.

Note that this recursive method reduces to Paschke's exact
sequence
(\ref{eq:Paschke})
in the case $n=2$.
When $n=2$, we are left with the single exact sequence
\begin{equation}
\begin{CD}
K_0(B_1) @>>> K_0(A \ltimes_\alpha {\mathbb Z}_n) @>>> K_0(A_0/I_0) \\
@AAA & & @VVV\\
K_1(A_0/I_0)  @<<<  K_1( A \ltimes_\alpha {\mathbb Z}_n) @<<< K_1(B_1)
\end{CD},
\end{equation}
where $B_1=B_{n-1}=A_0$ and $I_0=A_1^2$.
(Note that $\xi = -1$ in (\ref{eq:Asubk}), when $n=2$.)

\section{Example of $PSL_2(\mathbb R)$ acting on $H_2$}

Let $H_2$ be the hyperbolic plane modelled by the open upper half
place of ${\mathbb C}$.  Take $A=C_0(H_2)$.
The group $PSL_2({\mathbb R})$
acts on $H_2$ via
\begin{equation}
z \mapsto {ax + b \over cz + d}.
\end{equation}
The matrix
\begin{equation}
\begin{pmatrix}
0 &  1 \\ -1 & 1
\end{pmatrix}
\end{equation}
represents an element of order 3 in $PSL_2({\mathbb R})$.
We compute $K_*( A \ltimes_\alpha {\mathbb Z}_3)$,
where
\begin{equation}
{\mathbb Z}_3 = \bigl< \begin{pmatrix}
0 & 1 \\ -1 & 1 \end{pmatrix} \bigr>
\end{equation}
and
\begin{equation}
\bigl(\alpha f\bigr)(z)
= f\biggl(
\begin{pmatrix}
0 &  1 \\ -1 & 1
\end{pmatrix}^{-1}
z
\biggr)
= f\biggl(
\begin{pmatrix}
1 &  -1 \\ 1 & 0
\end{pmatrix}
z
\biggr)
= f(1-{1\over z}),
\end{equation}
 $f \in A$.

\vskip\baselineskip
\noindent{\bf Computation of $K_*(A_0)$.}\ 
Let
\begin{equation}
T= \begin{pmatrix}
1 & -1 \\ 1 & 0
\label{eq:exampleT}
\end{pmatrix}.
\end{equation}
Then the fixed-point algebra is
$A_0 = \{ \, f \in C_0(H_2) \, | \, f(Tz) = f(z) \, \}$.
\begin{figure}
\includegraphics[scale=0.6]{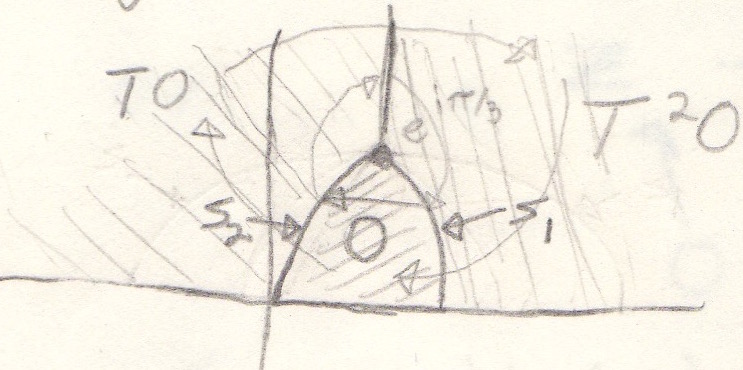} 
\caption{Drawing of region in the hyperbolic plane.
$O$ is a fundamental
region [Leh, 1966] for the action of T (\ref{eq:exampleT}).
}
\end{figure}
Find an open set $O \subseteq H_2$ (actually a polygon)
such that $O$, $TO$, and $T^2O$ are pairwise disjoint
and ${\overline {O \cup TO \cup T^2O}}=H_2$
(see Figure 1 and [Leh, 1966]).
The two sides $S_1$, $S_2$ of $O$ intersect in $T$'s only
fixed point $e^{i\pi/3}$.
$T$ maps each point in $S_1$ to the horizontally adjacent
point in $S_2$, and
$A_0 \cong \{\, f \in C_0(O\cup S_1\cup S_2)\,
| \, f(z) = f(Tz), z \in S_1\, \}.$
\begin{figure}
\includegraphics{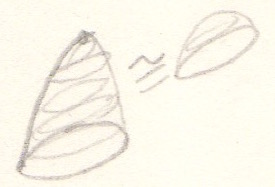} 
\caption{
Deformation of region in hyperbolic plane.
Identifying the sides
$S_1$ and $S_2$ folds the region $O$ into an upper hemisphere.
}
\end{figure}
Then $A_0$ is isomorphic to the algebra of continuous functions on the
upper half of the sphere $S^2$, vanishing at the equator (Figure 2).
Hence $A_0 \cong C_0({\mathbb R}^2)$ and
\begin{equation}
K_*(A_0) =
\begin{cases}
{\mathbb Z} & *=0 \\
0 & *=1. 
\end{cases}
\end{equation}

\vskip\baselineskip
\noindent{\bf Computation of $K_*(A_0/I_0)$.}\ 
Note that since $A$ is commutative, $A_1 A_2 = A_2 A_1$ and
therefore $I_0 = A_1 A_2 + A_2 A_1 = I_1 = A_1 A_2$.

Recall $A_i = \{\, f \in C_0(H_2)\,
| \, f(Tz) = \xi^i f(z) \}, i=0,1,2$ .
If $i \ne 0$,
every element of $A_i$ must vanish at any fixed point
of $T$.  Thus since functions in $A_0$ do not necessarily
vanish at $e^{i\pi/3}$ (recall that $e^{i \pi /3}$
is $T$'s only fixed point), $A_1 A_2 \subsetneq A_0$,
falling short of equality at the one point.
We have $A_0/I_0 \cong {\mathbb C}$, and therefore
\begin{equation}
K_*(A_0/I_0) =
\begin{cases}
{\mathbb Z} & *=0 \\
0 & *=1. 
\end{cases}
\end{equation}

\vskip\baselineskip
\noindent{\bf Use exact sequences to find 
$K_*(A \ltimes_\alpha {\mathbb Z}_3)$.}\ 
Now follow the recursive procedure
(\ref{eq:fullCrossCD})-(\ref{eq:nm2Exact}) outlined
above to find the $K$-theory of $A \ltimes_\alpha {\mathbb Z}_3$.
Our first exact sequence for $K_*(B_1) = K_*\biggl(\begin{pmatrix}
A_0 & A_1 \\ A_2 & A_0 \end{pmatrix}\biggr)$ is
\begin{equation}
\begin{CD}
K_0(A_0) @>>> K_0(B_1) @>>> K_0(A_0/I_0) \\
@AAA & & @VVV\\
K_1(A_0/I_0)  @<<<  K_1(B_1) @<<< K_1(A_0)
\end{CD}
\end{equation}
or
\begin{equation}
\begin{CD}
{\mathbb Z} @>>> K_0(B_1) @>>> {\mathbb Z} \\
@AAA & & @VVV\\
0  @<<<  K_1(B_1) @<<< 0
\end{CD}.
\end{equation}
So $K_0(B_1) = {\mathbb Z} \oplus {\mathbb Z}$, $K_1(B_1)=0$.
Our second exact sequence is
\begin{equation}
\begin{CD}
K_0(B_1) @>>> K_0(A \ltimes_\alpha {\mathbb Z}_3) @>>> K_0(A_0/I_0) \\
@AAA & & @VVV\\
K_1(A_0/I_0)  @<<<  K_1(A \ltimes_\alpha {\mathbb Z}_3) @<<< K_1(B_1)
\end{CD}
\end{equation}
or
\begin{equation}
\begin{CD}
{\mathbb Z} \oplus {\mathbb Z} @>>> K_0(A \ltimes_\alpha {\mathbb Z}_3)
@>>> {\mathbb Z} \\
@AAA & & @VVV\\
0  @<<<  K_1(A \ltimes_\alpha {\mathbb Z}_3) @<<< 0
\end{CD}.
\end{equation}
Thus $K_0(A \ltimes_\alpha {\mathbb Z}_3)= {\mathbb Z}^3$,
$K_1(A \ltimes_\alpha {\mathbb Z}_3)= 0$.

\vskip\baselineskip
\vskip\baselineskip
\section{References}
\smallskip
\smallskip
\footnotesize

\noindent[{\bf Bro, 1977}] 
\, L.G. Brown, 
{\it Stable isomorphism of hereditary subalgebras
of $C^\star$-algebras},
Pacific J. Math. {\bf 71(2)} (1977), 335-348.

\noindent[{\bf Ich, 1989}] 
\, R. Ichihara, 
{\it An extension of Paschke's theorem and
the Mayer-Vietoris type exact sequences},
Math. Japonica {\bf 34(1)} (1989), 27-33.

\noindent[{\bf Leh, 1966}]
\, J. Lehner,
{\it A short course in automorphic functions}, 
Holt, Rinehart, and Winston, New York, 1966.

\noindent[{\bf Pas, 1985}]
\, W.L. Paschke, 
{\it ${\mathbb Z}_2$-equivariant $K$-theory},
Lecture Notes in Math. {\bf 1132} (1985), 362-373.

\vskip\baselineskip

\noindent{Web Page: \url{http://www.svpal.org/~lsch/Math/indexMath.html}.}

\end{document}